\documentclass{amsart}% formato articolo bello
\usepackage[applemac]{inputenc}
\usepackage[T1]{fontenc}
\usepackage{amsfonts}
\usepackage{amsmath}
\usepackage{amssymb}
\usepackage{amsthm}
\usepackage{graphicx}
\usepackage{mathrsfs}
\usepackage{a4wide}%formato largo (margini)
\usepackage{dcpic,pictex} % serve per fare i grafici che commutano
\usepackage{hyperref} %serve per fare i riferimenti ipertestuali nel pdf
\usepackage[Sonny]{fncychap} %rende i titoli dei capitoli fichi

\DeclareMathAlphabet{\mathpzc}{OT1}{pzc}{m}{it}
\newcommand{\hh}{\mathbb{H}}
\newcommand{\hhc}{\mathbb{H}_\mathbb{C}}
\newenvironment{dimo}{\textit{Proof.\ }}{\begin{flushright}$\square$\end{flushright}}

\newtheorem{theorem}{Theorem}%[chapter]
\newtheorem{proposition}[theorem]{Proposition}
\newtheorem{corollary}[theorem]{Corollary}
\newtheorem{lemma}[theorem]{Lemma}
\theoremstyle{definition}
\newtheorem{defi}{Definition}%[chapter]
\newtheorem{exe}{Example}%[chapter]
\newtheorem{remark}{Remark}%[chapter]
\begin{document}
\title[Some properties for quaternionic slice-regular functions...]{Some properties for quaternionic slice-regular functions\\ 
on domains without real points}
\author{Amedeo Altavilla} 
\address{Department of Mathematics, University of Trento, Via Sommarive 14 I-38123 Italy}
\email{amedeo.altavilla@unitn.it}
\thanks{Department of Mathematics, University of Trento, Via Sommarive 14 I-38123 Italy}
\date{June 18, 2013}
\maketitle

\begin{abstract}
The theory of slice regular functions over the quaternions, introduced by Gentili and Struppa 
in \cite{gentilistruppa}, was born on domains that intersect the real axis. 
 This hypothesis can be overcome using the theory of stem functions introduced by Ghiloni and Perotti (\cite{ghiloniperotti}), in the context of real alternative algebras.
 In this paper I will recall the notion and the main properties of stem functions.
 After that I will introduce the class of slice regular functions induced by stem functions and, in this set, 
 I will extend the identity principle, the maximum and minimum modulus principles and the open mapping theorem. 
 Differences will be shown between the case when the domain does or does not intersect the real axis.\medskip

% \begin{keywords} function of a hypercomplex variable; slice regular functions; identity principle; minimum modulus principle; open mapping theorem\end{keywords}\medskip
% 
% \begin{classcode}30G35; 30C80 \end{classcode}\medskip

\end{abstract}

\section{Introduction}

Let $\mathbb{H}$ denote the algebra of quaternions. An element $x$ of $\mathbb{H}$ is of the form
$x=x_{0}+ix_{1}+jx_{2}+kx_{3}$, where $x_{l}$ are real numbers and $i,j,k$ are such that:
\begin{equation*}
i^{2}=j^{2}=k^{2}=-1,\quad ij=-ji=k, \, jk=-kj=-1, \,ki=-ik=j. 
\end{equation*}

An interesting and promising theory, based on a definition of regularity for quaternionic-valued functions of one
quaternionic variable given by C. G. Cullen in 1965 (see \cite{cullen}), was reintroducted and developed, in the last years, by G. Gentili, D. C. Struppa and others 
(see \cite{gentilistruppa}, \cite{genstostru} and their bibliography). More precisely the main concept is the following. Denote by $\mathbb{S}$ the sphere
of imaginary units:
\begin{equation*}
\mathbb{S}:=\{x\in\mathbb{H}\,|\,x^2=-1\},
\end{equation*}
 then a point $x$ in $\mathbb{H}$ can be written as $x=\alpha+\beta I$, where $\alpha,\beta\in\mathbb{R}$ and
 $I\in\mathbb{S}$. Therefore, putting $\mathbb{C}_{I}$ the real subspace of $\mathbb{H}$ generated by $1$ and
 $I$, we give the following definition:
 \begin{defi}
 Let $\Omega$ be a domain in $\mathbb{H}$ and let $f:\Omega\rightarrow \mathbb{H}$ be a real differentiable 
 function. $f$ is said to be Cullen-regular (briefly, regular), if, for all $I\in\mathbb{S}$, the function 
 $\overline{\partial} f:\Omega\cap\mathbb{C}_{I}\rightarrow \mathbb{H}$ defined by
 \begin{equation}\label{cullen}
 \overline{\partial} f(\alpha+\beta I)=\frac{1}{2}\left(\frac{\partial}{\partial	\alpha}+I\frac{\partial}{\partial\beta}\right)f|_{I}(\alpha+\beta I)
 \end{equation}
vanishes identically.
 \end{defi}
 Examples of such functions are convergent quaternionic power series defined on a ball centered in the origin
 \begin{equation*}
 f:B(0,R)\rightarrow \mathbb{H},\quad f(q)=\sum_{n\in\mathbb{N}}q^{n}a_{n},
 \end{equation*}
where $R>0$ denote the radius of convergence of the sum. It is also true the following theorem:
\begin{theorem}
Let $f:B(0,R)\rightarrow \mathbb{H}$ be a regular function, then there exists a sequence of quaternions
$\{a_{n}\}_{n\in\mathbb{N}}$ such that
\begin{equation*}
f(q)=\sum_{n\in\mathbb{N}}q^{n}a_{n}
\end{equation*}
for all $q\in B(0,R)$. In particular, $f\in \mathcal{C}^{\infty}(B(0,R))$.
\end{theorem}

 A beautiful description of this theory can be found in the recent 
monograph by G. Gentili, C. Stoppato and D. C. Struppa \cite{genstostru}, where, in particular are state the following results that can be found originally in \cite{gentilistruppa} (theorems \ref{teo1} and \ref{teo2}) and in \cite{gentilistoppato} 
(theorems \ref{teo3} and \ref{teo4}). But first we need a couple of definitions.

\begin{defi}
Let $\Omega$ be a domain in $\mathbb{H}$ that intersects the real axis. $\Omega$ is called a slice domain if,
for all $I\in\mathbb{S}$, the intersection $\Omega_{I}=\Omega\cap\mathbb{C}_{I}$ is a domain of 
$\mathbb{C}_{I}$.

Furthermore, a set $T\subset \mathbb{H}$ is called (axially) symmetric if, for all points $\alpha +\beta I\in T$,
the set $T$ contains the whole sphere $\alpha +\beta \mathbb{S}$.
\end{defi}

\begin{defi}
Let $f$ be a regular function on a symmetric slice domain $\Omega$. We define the degenerate set of 
$f$ as the union $D_{f}$ of the $2$-spheres $S=\alpha+\beta\mathbb{S}$ (with $\beta\neq 0$), such that $f|_{S}$
is constant.
\end{defi}

The following theorems hold.

\begin{theorem}\label{teo1}\textbf{\upshape(Identity Principle)}
Let $f$ be a regular function on a slice domain $\Omega$. If, for some $I\in\mathbb{S}$, $f$ is equal to
zero on a subset of $\Omega_{I}$ having an accumulation point in $\Omega_{I}$, then $f\equiv 0$ in $\Omega$.
\end{theorem}

\begin{theorem}\label{teo2}\textbf{\upshape(Maximum Modulus Principle)}
Let $\Omega\subset\mathbb{H}$ be a slice domain and let $f:\Omega\rightarrow \mathbb{H}$ be regular.
If $|f|$ has a relative maximum at $p\in\Omega$ then $f$ is constant. 
\end{theorem}

\begin{theorem}\label{teo3}\textbf{\upshape(Minimum Modulus Principle)}
Let $\Omega\subset\mathbb{H}$ be a slice domain and let $f:\Omega\rightarrow \mathbb{H}$ be regular.
If $|f|$ has a relative minimum at $p\in\Omega$ then either $f$ is constant or $f( p)=0$.
\end{theorem}

\begin{theorem}\label{teo4}\textbf{\upshape(Open Mapping Theorem)}
Let $f$ be a regular function on a symmetric slice domain $\Omega$ and let $D_{f}$  be its degenerate set. Then
$f:\Omega\setminus D_{f}\rightarrow \mathbb{H}$ is open.
\end{theorem}

The aim of this paper is to extend these results to the case in which the domain of definition of the function
does not intersect the real line. To do this we need the tools introduced by R. Ghiloni and A. Perotti in \cite{ghiloniperotti}\textbf{, and developed in other works such as \cite{ghiloniperotti2},} in  the
more general context of the real alternative algebras. The main instrument of this theory is the notion of stem function i.e.: a 
complex intrinsec function from a domain in $\mathbb{C}$ to the complexified algebra $\mathbb{H}\otimes_{\mathbb{R}} \mathbb{C}$.
With the notion of stem function it is possible to construct a reasonable class of quaternionic-valued functions of one
quaternionic variable defined on domains that could not intersects the real line: the set of slice functions. Then, in this class,
it is possible to extract a subset of regular functions that coincide to the set of Cullen-regular functions in the case in which the 
domain of definition is slice. More details about this constructions will be given in section 2. 

In section 3 we will extends the identity principle and define the set of slice-constant functions.

In section 4 we will give the extensions of the maximum and minimum modulus principle.

Finally in section 5 will be given the statement of the open mapping theorem for regular functions that are not defined on real points.

In sections 3,4 and 5 will be given examples to underline the differences between the case in which
the domain of definition of the function does or does not intersects the real line.

Let's begin then, with some introductory materials.

\section{Preliminaries}

Given an element $x=x_{0}+ix_{1}+jx_{2}+kx_{3}$ in the algebra of quaternions, we define its conjugate as $x^c=x_0-ix_1-jx_2-kx_3$.
It is clear then that
\begin{itemize}
 \item $(x^c)^c=x$;
 \item $(xy)^c=y^cx^c$;
 \item $x^c=x$, $\forall x\in\mathbb{R}$.
\end{itemize}
We remember also that, for every $x\in\mathbb{H}$, is defined its (squared) norm  as $n(x)=xx^{c}$.
Let now $\mathbb{H}_\mathbb{C}=\mathbb{H}\otimes_\mathbb{R}\mathbb{C}$ be the real tensor product between $\mathbb{H}$ and the complex plane. An element of $\mathbb{H}_{\mathbb{C}}$ is a sum 
$w=x+\sqrt{-1}y$, where $x,y\in\mathbb{H}$. In $\mathbb{H}_{\mathbb{C}}$ are then defined two conjugations:
\begin{itemize}
 \item $w^c=(x+\sqrt{-1}y)^c=x^c+\sqrt{-1}y^c$;
 \item $\overline{w}=\overline{x+\sqrt{-1}y}=x-\sqrt{-1}y$.
\end{itemize}

\begin{defi}
 Let $D\subset \mathbb{C}$ be an open set. If a function $F:D\rightarrow \mathbb{H}_{\mathbb{C}}$ is complex intrinsec, i.e. satisfies the condition
 \begin{equation}\label{intrinsec}
  F(\overline{z})=\overline{F(z)},\,\forall z\in D\, s.t.\,\overline{z}\in D,
 \end{equation}
then $F$ is called a \underline{stem function} on D.
\end{defi}

Due to the intrinsic behaviour of stem functions, there are no restriction to assume that $D$ is symmetric with respect to the real axis, i.e.:
\begin{equation*}
 D=conj(D):=\{z\in\mathbb{C}\,|\,\overline{z}\in D\}.
\end{equation*}
In fact, if this is not the case, $F$ can be extended to $D\cup conj(D)$ by imposing (\ref{intrinsec}).

\begin{remark}\label{evenodd}
 $F$ is a stem function if and only if the $\mathbb{H}$-valued components $F_1$, $F_2$ of $F=F_1+\sqrt{-1}F_2$ form an even-odd pair with respect to the imaginary part of $z$, i.e.:
 \begin{equation*}
  F_1(\overline{z})=F_1(z),\,F_2(\overline{z})=-F_2(z),\,\forall z\in D
 \end{equation*}
\end{remark}

\begin{remark}\label{cplxcurve}
In remark 3 of \cite{ghiloniperotti} is described, in a more general context the following construction:
as a real vector space $\mathbb{H}$ has dimension $4$, so, let $\mathcal{B}=\{u_{k}\}_{k=1}^{4}$ be a basis for $\mathbb{H}$. The function $F$ can be identified
with a complex intrinsic curve in $\mathbb{C}^{4}$. Let $F(z)=F_{1}(z)+\sqrt{-1}F_{2}(z)=\sum_{k=1}^{4}F_{\mathcal{B}}^{k}(z)u_k$, with $F_{\mathcal{B}}^{k}(z)\in \mathbb{C}$. Then 
\begin{equation*}
\tilde{F}_{\mathcal{B}}=(F_{\mathcal{B}}^{1},F_{\mathcal{B}}^{2},F_{\mathcal{B}}^{3},F_{\mathcal{B}}^{4}):D\rightarrow \mathbb{C}^{4}
\end{equation*}
satisfies $\tilde{F_{\mathcal{B}}}(\overline{z})=\overline{\tilde{F_{\mathcal{B}}}(z)}$. 
Giving to $\mathbb{H}$ the unique manifold structure as a real vector space, 
we get that a stem function $F$ is of class $C^{k}$ or real-analytic if and only if the same property holds for $\tilde{F_{\mathcal{B}}}$. 
Moreover this notion is independent of the choice of the basis of $\mathbb{H}$.\
\end{remark}

\begin{defi}
Given an open subset $D$ of $\mathbb{C}$ we define 
\begin{equation*}
\Omega_{D}:=\{x=\alpha +\beta J\in \mathbb{C}_{J}\,|\, \alpha,\beta\in\mathbb{R},\alpha + i\beta\in D, J\in\mathbb{S}\}.
\end{equation*}
We call these type of set \underline{circular set}.
\end{defi}

From now on $D$ will always be an open set of $\mathbb{C}$ and $\Omega_{D}$ will be its associated circular set.
% \begin{lemma}
%  If $D$ is a connected subset of $\mathbb{C}$ then $\Omega_D$ is connected in $\mathbb{H}$.
% \end{lemma}
% \begin{dimo}
%  Since the product of connected sets is connected, then $D\times \mathbb{S}$ is connected. The thesis follows because the map:
%  \begin{equation*}
%   \begin{array}{rlll}
%    \imath: & D\times\mathbb{S} & \rightarrow & \Omega_D\\
%     & (\alpha+i\beta,J) & \mapsto & \alpha+\beta J
%   \end{array}
%  \end{equation*}
% is continuous.
% \end{dimo}

In the following we will use the notations
\begin{equation*}
 D_J:=\Omega_D\cap\mathbb{C}_J,\quad D_J^+:=\Omega_D\cap\mathbb{C}_J^+,
\end{equation*}
where $\mathbb{C}_J:=\{x=\alpha +\beta J\in\mathbb{H}\,|\,\alpha,\beta\in\mathbb{R}\}$ and $\mathbb{C}_J^+:=\{x=\alpha +\beta J\in\mathbb{H}\,|\,\alpha,\beta\in\mathbb{R}\, ,\,\beta\geq 0\}$.

We are now in position to define slice functions:

\begin{defi}
Any stem function $F:D\rightarrow \mathbb{H}_{\mathbb{C}}$ induces a \underline{(left) slice function}
\begin{equation*}
f=\mathcal{I}(F):\Omega_{D}\rightarrow \mathbb{H}
\end{equation*}
defined as follows. If $x=\alpha + \beta J\in D_{J}$, we set
\begin{equation*}
f(x):=F_{1}(z)+JF_{2}(z),\quad z=\alpha+\sqrt{-1}\beta.
\end{equation*}
The set of (left) slice functions will be denoted by
\begin{equation*}
\mathcal{S}(\Omega_{D}):=\{f:\Omega_{D}\rightarrow \mathbb{H}\,|\, f=\mathcal{I}(F), F:D\rightarrow \mathbb{H}_{\mathbb{C}} \mbox{ stem function}\}.
\end{equation*}
\end{defi}
\textbf{The slice function $f$ is well defined, since $(F_1 , F_2 )$ is an even-odd pair w.r.t. $\beta$ and then
$f(\alpha+(-\beta)(-J))=F_1(\overline{z})+(-J)F_2(\overline{z})=F_1(z)+JF_2(z)$.}
Of course, there is an analogous definition for right slice functions when the element $J\in\mathbb{S}$ is placed on the right of $F_{2}(z)$.
\begin{remark}
 Since $\mathcal{I}(Fa+G)=\mathcal{I}(F)a+\mathcal{I}(G)$, for every $a\in\mathbb{R}$ and for every $F,G$ stem functions, 
 then $\mathcal{S}(\Omega_D)$ results to be a real vector space
\end{remark}
\textbf{Using stem functions to generate slice functions will be very useful in the next, since some 
computations wich are allowed in the set $\mathbb{H}_{\mathbb{C}}$ of complexified quaternions, are not in $\mathbb{H}$.}
\begin{exe}\label{sliceexe}
\begin{enumerate}
 \item Clearly the functions $z=Re(z)+\sqrt{-1}Im(z)$ and $\overline{z}=Re(z)-\sqrt{-1}Im(z)$ induces the functions $x$ and $x^c$ respectively.
 \item For any $a\in\mathbb{H}$ ,$F(z):= z^na = Re(z^n)a+\sqrt{-1}(Im(z^n)a)$ induces the monomial $f(x)=x^na\in\mathcal{S}(\mathbb{H})$.
 \item By linearity, we get all the standard polynomials $p(x) =\sum_{j=0}^n x^j a_j$ with right quaternionic coefficients. 
 More generally, every convergent power series $\sum_jx^ja_j$, with (possibly infinite) convergence radius $R$ (with respect to $|x|^{2} = n(x)$),
 belongs to the space $\mathcal{S}(B_R)$, where $B_R$ is the open ball of $\mathbb{H}$ centered in the origin with radius $R$.
 \item The two functions $G(z):=Re(z^n)a$ and $H(z):=\sqrt{-1}Im(z^n)a$ are complex intrinsic on $\mathbb{C}$.
 They induce respectively the slice functions $g(x)=Re(x^n)a$ and $h(x)=f(x)-g(x)=(x^n-Re(x^n))a=Im(x^n)a$.
 %The difference $G(z)-H(z)=\overline{z}^na$ induces $g(x)-h(x)=(2Re(x^n)-x^n)a=\overline{x}^na$.
\end{enumerate}

\end{exe}

An important property of slice functions is that they can be recovered by their values on two semislices. More precisely, we have the following theorem:

\begin{theorem}\label{representation}
Let $J,K\in\mathbb{S}$ with $J\neq K$. Then every $f\in \mathcal{S}(\Omega_{D})$ is uniquely determined
by its values on $D_J^+$ and $D_{K}^{+}$. More precisely we have the following formula
\begin{equation*}
f(x)=(I-K)((J-K)^{-1}f(\alpha +\beta J))-(I-J)((J-K)^{-1}f(\alpha +\beta K))
\end{equation*}
for all $I\in\mathbb{S}$, for all $x=\alpha + \beta I\in D_{I}$.
\end{theorem}

\textbf{The proof of this theorem can be found firstly in \cite{colgensabstru} where was proved for quaternionic regular
functions and was used to show an extension result. After that, Ghiloni and Perotti in \cite{ghiloniperotti} proved the same
theorem for slice functions which are not, in general, regular.}

The theorem does not exclude the possibility $K=-J$, and in this situation 
\begin{equation*}
f(x)=\frac{1}{2}(f(\alpha +\beta J)+f(\alpha -\beta J))-\frac{I}{2}(J(f(\alpha +\beta J)-f(\alpha -\beta J))).
\end{equation*}
\textbf{Moreover, if $I=J$, we are no more in the hypothesis of the theorem but we have the trivial equality}
\begin{equation*}
f(x)=\frac{1}{2}(f(x)+f(x^{c}))+\frac{1}{2}(f(x)-f(x^{c})),
\end{equation*}
where clearly $\frac{1}{2}(f(x)+f(x^{c}))=F_{1}(z)$ and $\frac{1}{2}(f(x)-f(x^{c}))=JF_{2}(z)$.

\begin{defi}
 We define the %\underline{spherical value} of $f$ in $x\in\Omega_{D}$
% \begin{equation*}
% v_{s}f(x):=\frac{1}{2}(f(x)+f(x^{c})),
% \end{equation*}
% and the 
\underline{spherical derivative} of $f$ in $x\in\Omega_{D}\setminus \mathbb{R}$
\begin{equation*}
\partial_{s}f(x):=\frac{1}{2}Im(x)^{-1}(f(x)-f(x^{c}))
\end{equation*}

\end{defi}

\begin{remark}\label{spherical}
We have that 
% $v_{s}f=\mathcal{I}(F_{1}(z))$ on $\Omega_{D}$ and 
$\partial_{s}f=\mathcal{I}(\frac{F_{2}(z)}{Im(z)})$ on $\Omega_{D}\setminus\mathbb{R}$.
Obviously this  function is constant on every sphere $\mathbb{S}_{x}=\{y\in\mathbb{H}\,|\,y=\alpha +\beta I,I\in\mathbb{S}\}$. In other terms:
\begin{equation*}
\partial_{s}(\partial_{s}(f))=0,
\end{equation*}
moreover $\partial_{s}f=0$ if and only if $f$ is constant on $\mathbb{S}_{x}$. 
%and, in this case, $f=v_{s}f(x)$ on $\mathbb{S}_{x}$.
If $\Omega_{D}\cap \mathbb{R}\neq\emptyset$, under some regularity hypotesis on $F$ (e.g.: differentiability of $F_2$), $\partial_{s}f$ can be 
extended continuously as a slice function on $\Omega_{D}$. 
% By definition the following identity holds for all $x\in\Omega_{D}$
% \begin{equation*}
% f(x)=v_{s}f(x)+Im(x)\partial_{s}f(x).
% \end{equation*}
\end{remark}

% There is a regularity result for slice functions depending on their stem functions
% 
% \begin{theorem}
% Let $f=\mathcal{I}(F)\in\mathcal{S}(\Omega_{D})$.
% \begin{enumerate}
% \item If $F\in C^{0}(D)$ then $f,v_{s}f\in C^{0}(\Omega_{D})$, $\partial_{s}f\in C^{0}(\Omega_{D}\setminus\mathbb{R})$;
% \item If $F\in C^{2s+1}(D)$, $s\in\mathbb{N}$, then $f,v_{s}f,\partial_{s}f\in C^{s}(\Omega_{D})$;
% \item If $F\in C^{\omega}$ then $f,v_{s}f,\partial_{s}f\in C^{\omega}(\Omega_{D})$.
% \end{enumerate}
% \end{theorem}

We will denote by $\mathcal{S}^{1}(\Omega_{D})$ the set of slice function induced by a $C^{1}$ stem function:
\begin{equation*}
\mathcal{S}^{1}(\Omega_{D}):=\{f=\mathcal{I}(F)\in\mathcal{S}(\Omega_{D})\,|\,F\in C^{1}(D)\}.
\end{equation*}

Let now $f\in\mathcal{I}(F)\in\mathcal{S}^{1}(\Omega_{D})$ and $z\in\alpha +\sqrt{-1}\beta\in D$, 
then $\frac{\partial F}{\partial \alpha}, \sqrt{-1}\frac{\partial F}{\partial \beta}$ are $C^{0}$ 
stem functions on $D$, and so also $\frac{\partial F}{\partial z}=\frac{1}{2}(\frac{\partial F}{\partial \alpha}-\sqrt{-1}\frac{\partial F}{\partial \beta})$ and
$\frac{\partial F}{\partial \overline{z}}=\frac{1}{2}(\frac{\partial F}{\partial \alpha}+\sqrt{-1}\frac{\partial F}{\partial \beta})$ are continuous.

\begin{defi}
Let $f=\mathcal{I}(F)\in\mathcal{S}^{1}(\Omega_{D})$. We define the following continuous slice functions
\begin{equation*}
\frac{\partial f}{\partial x}:=\mathcal{I}\left(\frac{\partial F}{\partial z}\right)
\end{equation*}
\begin{equation*}
\frac{\partial f}{\partial x^{c}}:=\mathcal{I}\left(\frac{\partial F}{\partial \overline{z}}\right).
\end{equation*}
\end{defi}

Left multiplication by $\sqrt{-1}$ defines a complex structure on $\mathbb{H}_{\mathbb{C}}$ and, with respect to this structure, a $C^{1}$ function
\begin{equation*}
F=F_{1}+\sqrt{-1}F_{2}:D\rightarrow\mathbb{H}_{\mathbb{C}}
\end{equation*}
is holomorphic if and only if satisfy the Cauchy-Riemann equations
\begin{equation*}
\frac{\partial F_{1}}{\partial \alpha}=\frac{\partial F_{2}}{\partial\beta}\, ,\,\frac{\partial F_{2}}{\partial\beta}=-\frac{\partial F_{2}}{\partial \alpha}\, ,\, z=\alpha +i\beta \in D
\end{equation*}
or equivalently if
\begin{equation*}
\frac{\partial F}{\partial \overline{z}}=0.
\end{equation*}
This condition is equivalent to require that, for any basis $\mathcal{B}$, the complex curve $\tilde{F}_{\mathcal{B}}$ defined in remark \ref{cplxcurve} is holomorphic.

We are now in position to define slice regular functions.
\begin{defi}
A function $f\in\mathcal{S}^{1}(\Omega_{D})$ is \underline{(left) slice regular} if its stem function $F$ is holomorphic. The set of slice regular function will be denoted by
\begin{equation*}
\mathcal{SR}(\Omega_{D}):=\{f\in\mathcal{S}^{1}(\Omega_{D})\,|\,f=\mathcal{I}(F),F:D\rightarrow \mathbb{H}_{\mathbb{C}}\mbox{ holomorphic}\}.
\end{equation*}
\end{defi}

The set of slice regular functions is again a real vector space.

\begin{exe}
 The polynomials and power series in the previous example \ref{sliceexe} are non trivial slice regular functions.
\end{exe}

The previous example shows that the theory of Cullen-regularity and slice regularity coincides when the domain
of definition is a slice domain. \textbf{Moreover in \cite{ghiloniperotti2}, Ghiloni and Perotti, showed, in the more general 
context of real alternative algebras, that asking regularity (in the sense of Cullen)
for a quaternionic function does not imply sliceness if the domain does not intersects the real axis. So, it seems that, the right extension of this theory
for domains which does not have real points, must be the one of regular functions which are slice.}

\begin{remark}
A function $f\in\mathcal{S}^{1}(\Omega_{D})$ is slice regular if and only if $\frac{\partial f}{\partial x^{c}}\equiv 0$. Moreover, if $f$ is slice regular, then 
also $\frac{\partial f}{\partial x}$ is slice regular on $\Omega_{D}$.
\end{remark}

\begin{remark}
Since $\partial_{s}f$ is $\mathbb{H}$-valued, then it is slice regular only when is locally constant.
\end{remark}

We will introduce now a useful notation: let $f=\mathcal{I}(F):\Omega_D\rightarrow \mathbb{H}$ then we denote the restriction over a complex plane or a complex half-plane, 
respectively, as
\begin{equation*}
 f_J:=f|_{D_J}:D_J\rightarrow \mathbb{H},\quad f_J^+:=f|_{D_J^+}:D_J^+\rightarrow \mathbb{H}.
\end{equation*}

\begin{proposition}
Let $f=\mathcal{I}(F)\in\mathcal{S}^{1}(\Omega_{D})$, then $f$ belong to $\mathcal{SR}(\Omega_{D})$ if and only if the restriction 
$f_{J}^+$ is holomorphic for every $J\in\mathbb{S}$,
with respect to the complex structures on $D_{J}$ and $\mathbb{H}$ defined by left multiplication by $J$.
\end{proposition}

\begin{remark}
The proof of the proposition and the even-odd character of the pair $(F_{1},F_{2})$ shows that, in order to get slice regularity of $f=\mathcal{I}(F)$
it is sufficient to assume that two restrictions $f_{J}^+$, $f_{K}^+$ ($J\neq K$) are holomorphic on $D_{J}^{+}$ and 
$D_{K}^{+}$ respectively. The possibility $K=-J$ is not excluded.
\end{remark}

We want now to multiplicate slice regular functions.
In general, the pointwise product of slice functions is not, a slice function, so we need another notion of product.
The following, \textbf{introduced by Gentili and Struppa in \cite{gentilistruppa} and by Ghiloni and Perotti in \cite{ghiloniperotti} in the context of real alternative agebras}, 
is the notion that we will use. 
\begin{defi}
Let $f=\mathcal{I}(F)$, $g=\mathcal{I}(G)$ $\in\mathcal{S}(\Omega_{D})$ the \underline{(slice) product} of $f$ and $g$ is the slice function
\begin{equation*}
f\cdot g:=\mathcal{I}(FG)\in\mathcal{S}(\Omega_{D}).
\end{equation*}
\end{defi}
 
 Sometimes the slice product between $f$ and $g$ is denoted by $f*g$ (see \cite{gentilistruppa} or \cite{gentilistoppato}) and called regular product.
 \begin{proposition}
 If $f,g\in\mathcal{SR}(\Omega_{D})$ then $f\cdot g\in \mathcal{SR}(\Omega_{D})$
 \end{proposition}
% 
% For this product it holds a sort of Leibniz's rule:
% \begin{equation*}
% \partial_{s}(f\cdot g)=(\partial_{s}f)(v_{s}g)+(v_{s}f)(\partial_{s}g)
% \end{equation*}
 
\begin{defi}
The slice function $f=\mathcal{I}(F)$ is called \underline{real} if the $\mathbb{H}$-valued components $F_{1}$, $F_{2}$ are real valued.
% (if and only if $v_{s}f$ and $\partial_{s}f$ are real)
\end{defi}

\begin{proposition}
A slice function $f=\mathcal{I}(F)$ is a \underline{real slice function} if and only if, for all $J\in \mathbb{S}$, $f(D_J)\subset \mathbb{C}_{J}$.
\end{proposition}

\begin{lemma}
Let $f=\mathcal{I}(F),g=\mathcal{I}(G)\in\mathcal{SR}(\Omega_{D})$, with $f$ real, then the slice function $h:\Omega_{D}\rightarrow \mathbb{H}$, defined by

\begin{equation*}
 h:=\mathcal{I}(FG)
\end{equation*}
is such that
\begin{equation*}
h(x)=f(x)g(x),
\end{equation*}
and, moreover, $h$ belongs to $\mathcal{SR}(\Omega_{D})$.
\end{lemma}

\begin{dimo}
The proof of this lemma can be found in a more general context in remark 7 of \cite{ghiloniperotti}:
in general $(f\cdot g)(x)\neq f(x)g(x)$. If $x=\alpha +\beta J$ belongs to $D_{J}=\Omega_{D}\cap \mathbb{C}_{J}$ and $z=\alpha + i\beta$, then
\begin{equation*}
(f\cdot g)(x)=F_{1}(z)G_{1}(z)-F_{2}(z)G_{2}(z)+JF_{1}(z)G_{2}(z)+JF_{2}(z)G_{1}(z),
\end{equation*}
while
\begin{equation*}
f(x)g(x)=F_{1}(z)G_{1}(z)+JF_{2}(z)JG_{2}(z)+F_{1}(z)JG_{2}(z)+JF_{2}(z)G_{1}(z).
\end{equation*}
If the components $F_{1}$, $F_{2}$ of the first stem function $F$ are real-valued, then $(f\cdot g)(x)=f(x)g(x)$ for every $x\in\Omega_{D}$.
\end{dimo}

\begin{corollary}\label{correal}
Let $f=\mathcal{I}(F),g=\mathcal{I}(G)\in\mathcal{SR}(\Omega_{D})$, with $g$ real, then the slice function $h:\Omega_{D}\setminus V(g)\rightarrow \mathbb{H}$, defined by

\begin{equation*}
 h=\mathcal{I}(G^{-1}F)
\end{equation*}
is such that
\begin{equation*}
h(x)=\frac{1}{g(x)}f(x),
\end{equation*}
and, moreover, $h$ belongs to $\mathcal{SR}(\Omega_{D})$.
\end{corollary}

We can now \textbf{remember the following definitions. The following appear for the first time in the work by Colombo, Gentili, Sabadini and Struppa \cite{colgensabstru}, 
but can be found also in \cite{gentilistoppato} and \cite{gentilistoppato2}. Later was
generalized by Ghiloni and Perotti in the context of real alternative agebras \cite{ghiloniperotti}}.

\begin{defi}\label{symm}
Let $f=\mathcal{I}(F)\in\mathcal{S}(\Omega_{D})$, then also $F^{c}(z)=F(z)^{c}:=F_{1}(z)^{c}+iF_{2}(z)^{c}$ is a stem function.
We set
\begin{itemize}
\item $f^{c}:=\mathcal{I}(F^{c})\in\mathcal{S}(\Omega_{D})$;
\item $CN(F):=FF^{c}$;
\item $N(f):=f\cdot f^{c}=\mathcal{I}(CN(F))$ \underline{symmetrization} or \underline{normal function} of $f$.
\end{itemize}
\end{defi}

The symmetrization of $f$ is sometimes denoted by $f^s$.

\begin{remark}
We have that $(FG)^{c}=G^{c}F^{c}$, and so $(f\cdot g)^{c}=g^{c}\cdot f^{c}$, i.e.:
\begin{equation*}
N(f)=N(f)^{c},
\end{equation*}
while, in general, $N(f^{c})\neq N(f)$
\end{remark}

Let now spent a few words about the zero locus of slice functions. For more details see \cite{ghiloniperotti} \textbf{and \cite{gentilistoppato2}}. We will use the following notation
\begin{equation*}
V(f):=\{x\in\mathbb{H} \,|\, f(x)=0\}.
\end{equation*}

\begin{proposition}
Let $f\in\mathcal{S}(\Omega_{D})$. The restriction $f_{|_{\mathbb{S}_{x}}}$ is injective or constant for all $x\in\Omega_{D}$.
In particular, either $\mathbb{S}_{x}\subseteq V(f)$ or $\mathbb{S}_{x}\cap V(f)$ is at most a singleton.
\end{proposition}

\begin{theorem}\label{zerostruc}
Let $f=\mathcal{I}(F)\in\mathcal{S}(\Omega_{D})$. Let $x=\alpha + \beta J\in\Omega_{D}$, $z=\alpha +i\beta\in D$. 
One of the mutually exclusive statements holds:
\begin{enumerate}
\item $\mathbb{S}_{x}\cap V(f)=\emptyset$;
\item $\mathbb{S}_{x}\subseteq V(f)$ (in this case $x$ is called a real ($x\in \mathbb{R}$) or spherical ($x\notin \mathbb{R}$) zero of $f$);
\item $\mathbb{S}_{x}\cap V(f)$ consists of a single, non-real point (in this case $x$ is called an $\mathbb{S}$-isolated non-real zero of $f$).
\end{enumerate}
These three possibilities correspond, respectively to the following properties of $F(z)\in\mathbb{H}_{\mathbb{C}}$:
\begin{enumerate}
\item $CN(F)(z)=F(z)F(z)^{c}\neq 0$;
\item $F(z)=0$
\item $F(z)\neq 0$ and $CN(F)(z)=0$.
\end{enumerate}
\end{theorem}

% \begin{remark}
% An $\mathbb{S}$-isolated non-real zero $x$ of $f$ is given by the formula
% \begin{equation*}
% x=\alpha +\beta K\, ,\,K:=-\frac{F_{1}(z)F_{2}(z)^{c}}{n(F_{2}(z))}\, ,\, CN(F)(z)=0.
% \end{equation*}
% This formula can be rewritten in the following way
% \begin{equation*}
% x=Re(x)-v_{s}f(x)(\partial_{s}f(x))^{-1}.
% \end{equation*}
% \end{remark}

\begin{corollary}\label{zeros}
The following statements hold:
\begin{enumerate}
\item a real slice function has no $\mathbb{S}$-isolated non-real zeros;
\item for all $f\in\mathcal{S}(\Omega_{D})$
\begin{equation*}
V(N(f))=\bigcup_{x\in V(f)} \mathbb{S}_{x}.
\end{equation*}
\end{enumerate}
\end{corollary}

\begin{theorem}
Let $\Omega_{D}$ be connected. If $f$ is slice regular and $N(f)$ does not vanish identically, then
\begin{equation*}
\mathbb{C}_{J}\cap\bigcup_{x\in V(f)}\mathbb{S}_{x}
\end{equation*}
is closed and discrete in $D_{J}$ for all $J\in\mathbb{S}$.
If $\Omega_{D}\cap\mathbb{R}\neq\emptyset$, then $N(f)\equiv 0$ if and only if $f\equiv 0$.
\end{theorem}

\begin{proposition}
Let $f,g\in\mathcal{S}(\Omega_{D})$. Then $V(f)\subset V(f\cdot g).$
\end{proposition}

\begin{proposition}\label{propzeros}
Let $x\in\mathbb{H}$. If $f,g\in\mathcal{S}(\Omega_{D})$, then it holds:
\begin{equation*}
\bigcup_{x\in V(f\cdot g)} \mathbb{S}_{x}=\bigcup_{x\in V(f)\cap V(g)}\mathbb{S}_{x}
\end{equation*}

\end{proposition}

At this point, the preliminary notions are fixed. Let now pass to the new results.

\section{Identity principle}
In this section we will prove an analogous of the identity principle for slice regular functions. 
In \cite{gentilistruppa} the authors prove the statement for slice regular functions defined on a domain that intersects the real axis. 
More precisely their statement, in our language, is the following:

\begin{theorem} [\textit{Gentili - Struppa}]
Let $\Omega_D\subset \mathbb{H}$ be a connected domain such that $\Omega_D\cap\mathbb{R}\neq \emptyset$, and let $f:\Omega_{D}\rightarrow \mathbb{H}$ be a slice regular function.
If there exists $I\in\mathbb{S}$ such that $D_{I}\cap V(f)$ has an accumulation point, then $f\equiv 0$ on $\Omega_{D}$. 
\end{theorem}

Now, it is not possible to generalize this theorem to the case in which the domain $D$ does not intersect the real line. In fact, there is a counterexample.
\begin{exe}
Let $J\in\mathbb{H}$ be fixed. The slice regular function defined on $\mathbb{H}\setminus \mathbb{R}$ by
\begin{equation*}
f(x)=1-IJ,\quad x=\alpha +\beta I \in\mathbb{C}_{I}^{+}
\end{equation*}
is induced  by a locally constant stem function and its zero set $V(f)$ is the half plane $\mathbb{C}_{-J}^{+}\setminus\mathbb{R}$.
The function can be obtained by the representation formula in theorem \ref{representation} by choosing the constant values $2$ on 
$\mathbb{C}_{J}^{+}\setminus\mathbb{R}$ and $0$ on $\mathbb{C}_{-J}^{+}\setminus \mathbb{R}$.
Now, by definition, it is clear that this function does not satisfies the thesis of the theorem.
\end{exe}

This example inspires us the following definition:
\begin{defi}
Let $f=\mathcal{I}(F)\in\mathcal{S}(\Omega_{D})$. $f$ is called \underline{slice constant} if the stem function $F$ is locally constant.
\end{defi}

\begin{remark}
 It is important to notice that a slice constant function can have multiple values on a semislice. Actually this could happen if the domain is not connected.
\end{remark}

\begin{proposition}
Let $f\in\mathcal{S}(\Omega_{D})$ be a slice constant function, then $f$ is slice regular.
\end{proposition}

\begin{dimo}
The proof is trivial because of the nature of the stem function that generate a slice constant function.
\end{dimo}

A simple characterization is given by the following theorem:
\begin{theorem}
Let $f=\mathcal{I}(F)\in\mathcal{SR}(\Omega_{D})$. Then $f$ is slice constant if and only if
\begin{equation*}
\frac{\partial f}{\partial x}=\mathcal{I}\left(\frac{\partial F}{\partial z}\right)\equiv 0.
\end{equation*}
\end{theorem}

\begin{dimo}
Let $F$ be locally constant, then in a connected component of $D$, $F=a+ib$, with $a,b\in\mathbb{H}$. This entails obviously that $\frac{\partial F}{\partial z}=0$.
Vice versa, let $f=\mathcal{I}(F)\in\mathcal{SR}(\Omega_{D})$ such that $\frac{\partial F}{\partial z}\equiv 0$, then, recalling the intrinsic curve in
remark \ref{cplxcurve}, we have in a connected component of $D$, that 
\begin{equation*}
\frac{\partial F}{\partial z}=\frac{\partial}{\partial z}\sum F_{\mathcal{B}}^{k}(z)u_k=0
\end{equation*}
and so $F_{\mathcal{B}}^{k}=c_k\in\mathbb{C}$, for every $k$, and so also $F=\sum_{k=1}^{4}F_{\mathcal{B}}^{k}u_{k}=c'\in \mathbb{H}$.
\end{dimo}

\begin{remark}
The previous theorem tells that if we have a slice constant function $f\in\mathcal{SR}(\Omega_{D})$, then, given $J\in\mathbb{S}$, if $x\in D_{J}\setminus \mathbb{R}$
\begin{equation*}
f(x)=a+ Jb=a+\frac{Im(x)}{|Im(x)|}b,\quad a,b\in\mathbb{H}.
\end{equation*}
\end{remark}

It is now clear that if we want to obtain an identity principle we must control the set of slice constant functions. 
The following theorem clarifies the situation. 

\begin{theorem} (\textbf{Identity Principle}) Let $\Omega_D$ be a connected open set of $\mathbb{H}$. Given $f=\mathcal{I}(F):\Omega_{D}\rightarrow \mathbb{H} \in \mathcal{SR}(\Omega_{D})$, with
$V(f)=\{x\in\Omega_{D}\,|\,f(x)=0\}$ be its zero locus. If there exists $K\neq J\in\mathbb{S}$ such that both 
$D_{K}^{+}\cap V(f)$ and $D_{J}^{+}\cap V(f)$ contain accumulation points, then $f\equiv 0$ on $\Omega_{D}$.

\end{theorem}

\begin{dimo}
Let $x_{J}=\alpha_{1}+\beta_{1}J$, $x_{K}$ be accumulation points on $D_{J}^+$ and $D_{K}^+$ respectively.
After having fixed a basis ${u_k}$ for $\mathbb{H}$ and putting $z=\alpha_1+\sqrt{-1}\beta_1$ we have that\footnote{see remark \ref{cplxcurve}}
\begin{equation*}
\begin{array}{rcl}
0=f(x_{J}) & = & F_{1}(\alpha_{1}+i\beta_{1})+JF_{2}(\alpha_{1}+i\beta_{1})\\
& = & \sum\limits_{k=1}^{4}F_{\mathcal{B}}^{k}(z)u_k,
\end{array}
\end{equation*}
and so  all the four components $F_{\mathcal{B}}^{k}$ are equal to zero, and since these are holomorphic, for the identity principle in the complex case, 
they are identically zero on $D_J^+$.
Replacing $K$ for $J$ in the previous formula, we obtain that $f$ is identically zero also on $D_{K}^+$.
We now obtain the thesis thanks to the representation formula in theorem \ref{representation}.
\end{dimo}

\section{Maximum and minimum modulus principles}

In this section we will generalize the maximum modulus principle stated in \cite{gentilistruppa}, 
to the case of regular function defined over circular domains that does not intersect (in general) the real line. 
Before this we need a lemma. The proofs of the lemma and the theorem, follow the argument in \cite{gentilistruppa}, with the corrections needed in our context.

\begin{lemma}\label{meanvalue}
If $f=\mathcal{I}(F)$ is a slice regular function, and if $I\in\mathbb{S}$, then $f_{I}^{+}$ has the mean value property. 
\end{lemma}

\begin{dimo}
If $x=\alpha+\beta I\in D_I^+$ then we know that $f(x)=F_{1}(z)+IF_{2}(z)$, with $z=\alpha+\sqrt{-1}\beta$. But then, for every
point $a=\gamma +\delta I\in D_I^+$, and all positive real number $r$ such that $\overline{B(a,r)}\subset D_I^+$ we have,
\begin{equation*}
\begin{array}{rcl}
\displaystyle\frac{1}{2\pi}\displaystyle\int_{0}^{2\pi}f(a+re^{I\vartheta})d\vartheta & = & \\
 & = & \displaystyle\frac{1}{2\pi}\displaystyle\int_{0}^{2\pi}(F_{1}(a'+re^{\sqrt{-1}\vartheta})+IF_{2}(a'+re^{\sqrt{-1}\vartheta}))d\vartheta \\
  & =& F_{1}(a')+IF_{2}(a')=f(a)
\end{array}
\end{equation*}
where $a'=\gamma+\sqrt{-1}\delta$, and the penultimate equality holds because, 
restricting to $D_I^+$ and passing through the complex curve in remark \ref{cplxcurve} $F_{1}(z)+\sqrt{-1}F_{2}(z)=\sum_{k=1}^{4}F_{\mathcal{B}}^{k}(z)u_k$, 
with $F_{\mathcal{B}}^{k}(z)\in \mathbb{C}$, for each component we are in the hypothesis of the mean value property in the holomorphic case.
\end{dimo}

\begin{theorem}(\textbf{Maximum Modulus Principle})
Let $f=\mathcal{I}(F)\in\mathcal{SR}(\Omega_{D})$ with $\Omega_D$ connected circular domain. If there exists $J\neq K\in\mathbb{S}$ such that $|f_J^+|$ has relative maximum in  
$a\in D_J^+$ and $|f_K^+|$ has relative maximum in $b\in D_K^+$, then $f$ is slice-constant on $\Omega_{D}$.
\end{theorem}

\begin{dimo}
If $f(a)=f(b)=0$ the result is trivial. 
We will assume that at least one between $f(a)$ and $f(b)$ is different from zero. 
Let then $f(a)\neq 0$, using Lemma \ref{meanvalue} and following the proof in \cite{gentilistruppa}, we get that $f$ is constant on $D_J^+$. 
Now, if $f(b)=0$, the result is again trivial thanks to the representation formula in \ref{representation}. If $f(b)\neq 0$, repeating the argument in \cite{gentilistruppa}, 
we get that $f$ is constant also over $\Omega_{D}\cap\mathbb{C}_{K}^{+}$. The proof is concluded thanks to the representation formula in \ref{representation}.
\end{dimo}

\begin{remark}\label{counterex}
It must be noticed that the hypothesis of double relative maximum on two different semislices of $\Omega_{D}$ is not removable. Indeed there is a counterexample:  
let $J\in\mathbb{S}$ be fixed then the function $f:(\mathbb{H}\cap B(0,r))\setminus \mathbb{R}\rightarrow\mathbb{H}$ defined by
\begin{equation*}
f(x)=x+c-(x-c)IJ,\quad x=\alpha+\beta I,\, \mathbb{R}\ni c>r
\end{equation*}
is constant and equal to $2c$ on $\mathbb{C}_{-J}^{+}\setminus\mathbb{R}$, but is equal to $2x$ on $\mathbb{C}_{J}^{+}\setminus\mathbb{R}$.
\end{remark}

We have the following trivial corollary

\begin{corollary}
Let $f=\mathcal{I}(F)\in\mathcal{SR}(\Omega_{D})$. If there exists $J\neq K\in\mathbb{S}$ such that $|f_J^+|$ has relative maximum in 
$a\in\Omega_{D}\cap\mathbb{C}_{J}^{+}$ and $|f_K^+|$ has relative maximum in $b\in\Omega_{D}\cap\mathbb{C}_{K}^{+}$ and $f(a)=f(b)$, then $f$ is constant on $\Omega_{D}$.
\end{corollary}

\begin{remark}
\textbf{ The previous formulations of the maximum modulus and the identity principles can be generalized to the context of real alternative aglebras.
 The proofs, indeed, does not use sthe particular quaternionic structures or properties. For the identity principle it is very clear since the existence 
 of a complex holomorphic curve is always guarantee for a slice function with values and defined in a real alternative algebra (see \cite{ghiloniperotti} for more details).
 The same must be true for the maximum modulus principle. Moreover we are not so sure that the following minimum modulus principle and open mapping theorem
 can be generalized as they are. We are pretty convinced that, adding some admissibility hypothesis there are chances to have an analogous minimum modulus principle. Furthermore, 
 other consideration are needed to obtain a reasonable generalization of the open mapping theorem but we will not spent any word about this since the present work
 is about quaternionic functions. We leave to the future the right considerations for solving these questions.}
\end{remark}

Our next goal now is to obtain a minimum modulus principle that generalize the one in \cite{gentilistoppato} to
the case of domains without real points. This will enable us to prove the open mapping theorem. 
For this reason we need some introductory material regarding the symmetrization and the reciprocal of a slice regular function.

\begin{theorem}
Let $f\in\mathcal{SR}(\Omega_{D})$. Let $x,y\in\mathbb{R}$, such that $x+y\mathbb{S}\subset \Omega_{D}$. The zeros of $f^{c}$ on
$x+y\mathbb{S}$ are in bijective correspondence with those of $f$. Moreover $N(f)$ vanishes exactly on the set $x+y\mathbb{S}$
on which $f$ has a zero.
\end{theorem}

The previous theorem can be found in \cite{gentilistoppato} in the case of domains with
real points. In general follows combining proposition \ref{propzeros} and corollary \ref{zeros}.

\textbf{For the next definition we give again the same references as for definition \ref{symm}, but we underline that, so far, there is not a generalization in
the set of real alternative algebra. The following, indeed, is also a proposal for that purpose.}
\begin{defi}
Let $f=\mathcal{I}(F)\in\mathcal{SR}(\Omega_{D})$. We call the \underline{(regular) reciprocal} of $f$ the slice function
\begin{equation*}
f^{-\bullet}:\Omega_{D}\setminus V(N(f))\rightarrow \mathbb{H}
\end{equation*}
defined by	
\begin{equation*}
f^{-\bullet}=\mathcal{I}(CN(F)^{-1}F^c)
\end{equation*}
\end{defi}

From the previous definition it follows that, if $x\in\Omega_D$, then

\begin{equation*} 
f^{-\bullet}(x)=(N(f)(x))^{-1}f^{c}(x).
\end{equation*}

The regularity of the reciprocal just defined follow thanks to corollary \ref{correal}.

\textbf{The following propositions already appeared in \cite{gentilistoppato}, but, since the hypothesis about the domains is quite different, we
propose new proofs.}

\begin{proposition}
Let $f\in\mathcal{SR}(\Omega_{D})$ such that $V(f)=\emptyset$, then $f^{-\bullet}\in\mathcal{SR}(\Omega_{D})$ and 
\begin{equation*}
f\cdot f^{-\bullet}=f^{-\bullet}\cdot f=1.
\end{equation*}
\end{proposition}

\begin{dimo}
Since $V(f)=\emptyset$ then $V(N(f))=\emptyset$. So $(N(f))^{-1}$ and $f^{-\bullet}$ are well defined and regular on the whole $\Omega_{D}$.
We may consider then their regular product with other regular functions $g:\Omega_{D} \rightarrow \mathbb{H}$.
For all $g$, $(N(f)(x))^{-1}g(x)=((N(f))^{-1}\cdot g)(x)$. Than we have
\begin{equation*}
f^{-\bullet}\cdot f=(N(f))^{-1}\cdot f^{c}\cdot f=(N(f))^{-1}N(f)=1
\end{equation*}
and
\begin{equation*}
f\cdot f^{-\bullet}=f\cdot(N(f))^{-1}\cdot f^{c}=(N(f))^{-1}\cdot f\cdot f^{c}=(N(f))^{-1}N(f)=1.
\end{equation*}
\end{dimo}

\begin{proposition}\label{prodcomp}
Let $f,g\in\mathcal{SR}(\Omega_{D})$ then, for any $x\in\Omega_{D}\setminus V(f)$
\begin{equation*}
(f\cdot g)(x)=f(x)g(f(x)^{-1}xf(x)).
\end{equation*}
\end{proposition}

\begin{dimo}
Let $x=\alpha +J\beta$, with $J\in\mathbb{S}$ and $z=\alpha + i\beta$. We have
\begin{eqnarray*}
(f\cdot g)(x) & = & \mathcal{I}(FG)(x) \\
& = & F_{1}(z)G_{1}(z)-F_{2}(z)G_{2}(z)+J(F_{1}(z)G_{2}(z)+F_{2}(z)G_{1}(z))\\
& = & (F_{1}(z)+JF_{2}(z))(G_{1}(z)+\\
& & +(F_{1}(z)+JF_{2}(z))^{-1}(JF_{1}(z)G_{2}(z)-F_{2}(z)G_{2}(z)))\\
& = & (F_{1}(z)+JF_{2}(z))(G_{1}(z)+\\
 & & (F_{1}(z)+JF_{2}(z))^{-1}J(F_{1}(z)+JF_{2}(z))G_{2}(z))
\end{eqnarray*}
but $J'=(F_{1}(z)+JF_{2}(z))^{-1}J(F_{1}(z)+JF_{2}(z))=f(x)^{-1}Jf(x)\in\mathbb{S}$. So if we call $x'=\alpha+J'\beta=f(x)^{-1}xf(x)$ we obtain the thesis
\begin{equation*}
(f\cdot g)(x)=f(x)g(x').
\end{equation*}
\end{dimo}

Thanks to the previous proposition we have the following corollary

\begin{corollary}
Let $f\in\mathcal{SR}(\Omega_{D})$, then if we set $T_{f}(x)=f^{c}(x)^{-1}xf^{c}(x)$, we have
\begin{equation*}
f^{-\bullet}(x)=f(T_{f}(x))^{-1},
\end{equation*}
for all $x\in\Omega_{D}\setminus V(f)$.
\end{corollary}

The proof of the corollary, that is identical to the one in \cite{gentilistoppato}, 
is a trivial application of the proposition \ref{prodcomp} to the function $f^{-\bullet}$ remembering that:
\begin{equation*}
f^{-\bullet}(x)=N(f)(x)^{-1}f^{c}(x)=((f^{c}\cdot f)(x))^{-1}f^{c}(x).
\end{equation*}

\begin{proposition}
Let $f\in\mathcal{SR}(\Omega_{D})$, then $T_{f}$ and $T_{f^{c}}$ are mutual inverses w.r.t. composition. 
Moreover $T_{f}:\Omega_{D}\setminus V(f^c)\rightarrow \Omega_{D}\setminus V(f)$ is a diffeomorphism.
\end{proposition}

Again, the proof of this statement is in \cite{gentilistoppato}.

We recall now the definition of the degenerate set of a function.

\begin{defi}
 Let $f\in\mathcal{S}(\Omega_D)$ and let $x,y\in\mathbb{R}$, $y>0$ be such that $S=x+y\mathbb{S}\subset \Omega_D$. The 2-sphere $S$ is said to be \underline{degenerate}
 for $f$ if the restriction $f|_S$ is constant. The union $D_f$ of all degenerate spheres for $f$ is called \underline{degenerate set} of $f$.
\end{defi}

\begin{proposition}
 Let $f$ a slice function over $\Omega_D$, then we have the following equality:
 \begin{equation*}
  D_f=Ker(\partial_s f).
 \end{equation*}
 Moreover, if $f$ is nonconstant $D_f$ is closed in $\Omega_{D}\setminus\mathbb{R}$.
\end{proposition}

\begin{dimo}
 The proof of the statement is trivial thanks to remark \ref{spherical}. 
\end{dimo}

\begin{proposition}
If $f\in\mathcal{SR}(\Omega_{D})$ is nonconstant, then the interior of $D_{f}$ is empty.
\end{proposition}
\begin{dimo}
If ad absurdum there exist a point $p\in D_{f}$ and a neigborhood $\Omega_{U}$ of $p$ such that 
$\Omega_{U}\subset D_{f}$ and $f$ is nonconstant, then, for all $\alpha +\beta J\in\Omega_{U}$
\begin{equation*}
f(\alpha+\beta J)=F_{1}(\alpha+i\beta).
\end{equation*}
Since $f$ is slice regular we have that
\begin{equation*}
0=\frac{\partial F_{1}}{\partial \overline{z}}=\frac{\partial F_{1}}{\partial \alpha}+i\frac{\partial F_{1}}{\partial\beta},
\end{equation*}
but then $\frac{\partial F_{1}}{\partial \alpha}=0$ and $\frac{\partial F_{1}}{\partial\beta}=0$ separately and so
$F_{1}$ is equal to a constant in all $\Omega_{U}\subset D_{f}\subset\Omega_{D}$.
 Thanks to the identity principle we obtain that $f$ is constant.
\end{dimo}

\begin{remark}
 If $f$ is a slice regular function defined on $\Omega_D$ and  $S=x+y\mathbb{S}\subset \Omega_D$ ($x,y\in\mathbb{R}$, $y>0$), is not degenerate, then the restriction 
 $f|_S$ is a nonconstant affine map of $S$ onto a 2-sphere $b+\mathbb{S}c$ with $b,c\in\mathbb{H}$.
\end{remark}

Thanks to this remark we have the following proposition:

\begin{proposition}
 Let $f$ be a slice regular function defined on $\Omega_D$ let $x,y\in\mathbb{R}$, $y>0$ be such that 
 $S=x+y\mathbb{S}\subset \Omega_D$, $S\nsubseteq D_f$. Then $|f|_S|$ has a global minimum, a 
 global maximum and no other extremal point.
\end{proposition}

We are now ready to state a formulation of the minimum modulus principle.

\begin{theorem}\label{minmod}(\textbf{Minimum Modulus Principle}) 
Let $\Omega_{D}$ be a connected circular domain and let $f:\Omega_D\rightarrow \mathbb{H}$ be a slice regular function. 
If $|f_I^+|$ has a local minimum point $p=x+yI\in D_I^+$ 
 then either $f(p)=0$ or exists a $J\in\mathbb{S}$ such that $f_J^+$ is constant.
\end{theorem}

\begin{dimo}
 Suppose $f$ does not have zeroes in $S=x+y\mathbb{S}$. The reciprocal $f^{-\bullet}$ is defined on $\Omega_D\setminus V(N(f))$ which includes $S$. 
 Since $|f^{-\bullet}(q)|=1/|f(T_f(q))|$ for every $q$ and $T_f$ is a diffeomorphism, the fact that $|f_I^+|$ has a minimum at $p=T_f(p')$ implies that $|f\circ T_f|$ has a minimum
 at $p'=x'+y'I'$. As a consequence, $|f^{-\bullet}|$ has a maximum at $p'$. Now, by the maximum modulus principle $f^{-\bullet}$ is constant on 
 $D_I^+$.
 
 Suppose now that, for every $J\in\mathbb{S}$, $f_J^+$ is non-constant and $|f_I^+|$ has a minimum at $p=x+yI$. There must exists a point
 $p'\in S$ such that $f(p')=0$ and then $|f|$ has a minimum at $p'$. By the previous proposition , $|f|$ cannot have two distinct local minimum points on the same sphere $S$,
 unless $S$ is degenerate for $f$. As a consequence, either $f$ is constant on $S$ or $p=p'$. In both cases, $f(p)=f(p')=0$.
\end{dimo}

It is clear that this theorem can be refined adding some hypothesis. For instance, if one ask for $f$ to have two minimal points $p\neq q$  for its modulus that are sent by $T_f$ on two different semislices,
then one can conclude that $f$ is slice constant. This case could happen for example when $p,q$ belongs to the same sphere, because we know that $T_f$ maps any 2-sphere to itself.
Anyway this formulation of the minimum modulus principle is sufficient to prove the open mapping theorem.

\section{Open mapping theorem}

Let $\Omega_D$ a connected circular domain of $\mathbb{H}$. Given $f\in\mathcal{SR}(\Omega_D)$ we want to describe the set $K$ of semislices where the slice derivative 
$\frac{\partial f}{\partial x}$ is equal to zero. 
% Thanks to the identity principle we obtain that the kernel of the slice derivative of $f$ is the union of two
% (eventually empty) subsets
% 
% \begin{equation*}
%  V\left(\frac{\partial f}{\partial x}\right)=K^0\cup K^2,
% \end{equation*}
% where $K^0$ is a discrete subset of $\Omega_D$ and $K^2$ is a subset of semislices of $\Omega_D$. 
In particular, $K$ can be empty, a single semislice or the entire 
$\Omega_D$. In fact if two different semislices belongs to $K$, for the identity principle, $\frac{\partial f}{\partial x}$ is identically zero.
We now want to characterize the subset $K$. To do this, fix $I\in\mathbb{S}$ and $f\in\mathcal{SR}(\Omega_D)$. Let define the following integral
\begin{equation*}
 \int_{D_I^+}|f|:=\int_{D^+}|F_1+IF_2|d\mu,
\end{equation*}
where $d\mu$
  denote the lebesgue measure over $D^+$. The last integral is equal to
zero if and only if $|F_1+IF_2|=0$ on $D^+$, i.e.: if $|f_I^+|=0$. Let define the following operators
\begin{eqnarray*}
 \mathcal{L}:\mathbb{S}\times \mathcal{SR}(\Omega_D) & \rightarrow & \mathbb{R}\\
 (I,f) & \mapsto & \int_{D_I^+}|f|,
\end{eqnarray*}
and
\begin{eqnarray*}
 \mathcal{L}_f:\mathbb{S} & \rightarrow & \mathbb{R}\\
 I & \mapsto & \int_{D_I^+}|f|,
\end{eqnarray*}
then, introducing the following notation 
\begin{equation*}
D_A^+:=\bigcup_{I\in A}D_I^+,\quad \forall A\subset\mathbb{S},
\end{equation*}
we finally have that
$$D_{Ker(\mathcal{L}_{\frac{\partial f}{\partial x}})}^{+}=K.$$ 
We can now state the open mapping theorem: 

\begin{theorem}(\textbf{Open Mapping Theorem}) Let $f:\Omega_D\rightarrow \mathbb{H}$ be a slice regular function. Then 
$$f:\Omega_D\setminus (Ker(\partial_s f)\cup D_{Ker(T_{\frac{\partial f}{\partial x}})}^+)\rightarrow \mathbb{H}$$ is open. 
\end{theorem}

The proof of this theorem follow the one in complex case.

\begin{dimo}
 Let $U$ be an open set of $\Omega_D\setminus (Ker(\partial_s f)\cup D_{Ker(T_{\frac{\partial f}{\partial x}})}^+)$, the thesis is that $f(U)$ is open in $\mathbb{H}$. Let $p_0\in f(U)$, then there exist $q_0\in U$ such that 
 $p_0=f(q_0)$. Clearly, the function $f(q)-p_0$ vanishes in $q_0$. Now, theorem \ref{zerostruc} tells that either $q_0$ is an isolated zero or is part of a sphere
 $S$ where the function vanishes identically. Since by hypotesis we have removed the degenerate set from the domain of the function, the last option cannot hold and
 $q_0$ is an isolated zero for $f$.
 We have then that there exists an open ball $B=B(q_0,r)$ such that $\overline{B}\subset U$ and $f(q)-p_0\neq 0$ for all $q\in\partial B$, i.e.: there exists $\epsilon >0$
 such that $|f(q)-p_0|\geq 3\epsilon$ for all $q\in\partial B$. We choose now an arbitrary $p$ such that $|p-f(q_0)|=|p-p_0|<\epsilon$ and we have the following inequality:
 \begin{equation*}
  |f(q)-p|\geq |f(q)-p_0|-|p-p_0|\geq3\epsilon-\epsilon=2\epsilon,\quad \forall q\in\partial B.
 \end{equation*}
We have obtained that the minimum of $|f(q)-p|$ in $\overline{B}$ is strictly less then its minimum in $\partial B$, and so $|f(q)-p|$ must have a minimum in $B$.
By theorem \ref{minmod}, either $f(q)-p$ vanishes at the point of minimum or there exists a semislice where the function is constant. Since, by hypothesis, $f$ is 
non-constant in every semislice, then there exists a point $q\in B\subset U$ such that $f(q)=p$ and $p\in f(U)$ and the proof is concluded.
\end{dimo}

The hypothesis for which $f$ cannot be constant in any semislice is not removable. Indeed the counterexample in remark \ref{counterex}
give an information also in this direction. Let's define $f$ as in remark \ref{counterex}:
\begin{equation*}
 f:\mathbb{H}\setminus \mathbb{R}\rightarrow \mathbb{H}
 \end{equation*}
 \begin{equation*}
 f(x)=x-xIJ,\quad x=\alpha +\beta I
 \end{equation*}
$f$ is non-constant in every semislice except for $\mathbb{C}_{-J}^{+}$ in which is identically equal to zero. Moreover in every semislice
except for $\mathbb{C}_{-J}^{+}$ it assumes non-purely real values. 
We will prove that $f:\mathbb{H}\setminus \mathbb{R}\rightarrow \mathbb{H}$ is not open while $f:\mathbb{H}\setminus \mathbb{C}_{-J}^+\rightarrow \mathbb{H}$ it is.

Given then $B=B(y,r)$, with 
$y\in(\mathbb{H}\setminus\mathbb{R})\cap\mathbb{C}_{-J}^{+}$, $r>0$ such that $B\cap\mathbb{R}=\emptyset$, we will prove that $f(B)$ is not open.
But this is trivial because $f(B)=\{0\}\sqcup D$, where $D\subset \mathbb{H}\setminus \mathbb{R}$ and the union is disjoint 
(this is true because, writing the function explicitely, it is easy to see that $f(\mathbb{H}\setminus \mathbb{C}_{-J}^+)\cap\mathbb{R}=\emptyset$). This entails 
that it is not possible to find a ball $B'\subset \mathbb{H}$ such that $f(y)\in B'\subset f(B)$ because, otherwise, the intersection 
$B'\cap\mathbb{R}$ must be equal to an interval $(0-\epsilon ,0+\epsilon)$, for some $\epsilon >0$ but $f(B)\cap\mathbb{R}=\{0\}$.

Let see now that the same function $f$ restricted to $\mathbb{H}\setminus (\mathbb{R}\cup\mathbb{C}_{-J}^+)$ is open.
First of all, if $x=\alpha +\beta I$, then
\begin{equation*}
 f(x)=\alpha(1+I\cdot J)+\beta I+\beta J-\sqrt{1-(I\cdot J)^2}I\wedge J,
\end{equation*}
where $I\cdot J$ and $I\wedge J$ denote the scalar and the vectorial products respectively in $\mathbb{R}^3$.
But then again,
\begin{equation*}
 f(x)=\alpha(1+I\cdot J)+(2\beta^2+\alpha^2(1-(I\cdot J)^2))\left(\frac{\beta I+\beta J-\sqrt{1-(I\cdot J)^2}I\wedge J}{2\beta^2+\alpha^2(1-(I\cdot J)^2)}\right).
\end{equation*}
Now, the sets of the form  
\begin{equation*}
 A=(\alpha-\epsilon,\alpha+\epsilon)+(\beta-\delta,\beta+\delta)B_I(R)\in\mathbb{H}\setminus (\mathbb{R}\cup\mathbb{C}_{-J}^+),
\end{equation*}
with $B_I(R)=B(I,R)\cap \mathbb{S}^2$, form a basis for the topology in $\mathbb{H}\setminus (\mathbb{R}\cup\mathbb{C}_{-J}^+)$.
So to prove that $f$ is open we need to prove that $f(A)$ is open. 
Is clear that if we stay far from $\mathbb{C}_J^+$, then $\alpha(1+I\cdot J)$ sends the set $A$ in an open intervall.
It is also clear that $(2\beta^2+\alpha^2(1-(I\cdot J)^2))$ sends $A$ into an open intervall since $\beta-\delta>0$.
For the last part, since for $J\notin B_I(R)$, $I$, $J$ and $I\wedge J$ are linear independent, then the function 
$\frac{\beta I+\beta J-\sqrt{1-(I\cdot J)^2}I\wedge J}{2\beta^2+\alpha^2(1-(I\cdot J)^2)}$ send $A$ into an open set.
If $I=J$ we have not problem, since the image of $\frac{\beta I+\beta J-\sqrt{1-(I\cdot J)^2}I\wedge J}{2\beta^2+\alpha^2(1-(I\cdot J)^2)}$ contain 
a ball centered in $J$.

We conclude this paper with the following remark.

\begin{remark}
 The last example shows also that, if the domain $\Omega_D$ of definition of a nonconstant slice regular function $f$ does not contains 
 real points, then, in general,
 the set $f(\Omega_D)$ is not open in $\mathbb{H}$.
\end{remark}


\begin{thebibliography}{10}

\bibitem{ahlfors}
Ahlfors, LV. Complex analysis. McGraw-Hill Book Co. New York. 1978.

\bibitem{colgensabstru}
Colombo F, Gentili G, Sabadini I, Struppa DC.
Extension results for slice regular functions of a quaternionic variable. (English summary)
Adv. Math. 222 (2009), no. 5, 1793–1808. 

\bibitem{cullen}
Cullen CG. An integral theorem for analytic intrinsic functions on quaternions. Duke Math J.. 1965; 32: pp. 139--148 

\bibitem{genstostru}
Gentili G, Stoppato C, Struppa DC. Regular Functions of a Quaternionic Variable. Springer London, Limited.
Springer Monographs in Mathematics 2013.

\bibitem{gentilistoppato}
Gentili G, Stoppato C. The open mapping theorem for regular quaternionic functions. Ann. Sc. Norm. Super. Pisa Cl. Sci.. 2009; 5: pp. 805--815.

\bibitem{gentilistoppato2}
Gentili G, Stoppato C. Zeros of regular functions and polynomials of a quaternionic variable. Michigan Math. J. 56 (2008), no. 3, 655–667. 

\bibitem{gentilistruppa}
Gentili G, Struppa DC. A new theory of regular functions of a quaternionic variable. Adv. Math.. 2007; 216: pp. 279--301.

\bibitem{ghiloniperotti}
Ghiloni R, Perotti A. Slice regular functions on real alternative algebras. Adv. Math.. 2011; 226: pp. 1662--1691.

\bibitem{ghiloniperotti2}
Ghiloni R, Perotti A. Global differential equations for slice regular functions. Math. Nachr.. 2013; in press.

\bibitem{rudin}
Rudin W. Real and complex analysis. McGraw-Hill Book Co. New York.1987.

\end{thebibliography}
\end{document}